# Student's $t$-test for scale mixture errors

## Gábor J. Székely[1]

*Bowling Green State University, Hungarian Academy of Sciences*

**Abstract:** Generalized t-tests are constructed under weaker than normal conditions. In the first part of this paper we assume only the *symmetry* (around zero) of the error distribution (i). In the second part we assume that the error distribution is a *Gaussian scale mixture* (ii). The optimal (smallest) critical values can be computed from generalizations of Student's cumulative distribution function (cdf), $t_n(x)$. The cdf's of the generalized $t$-test statistics are denoted by (i) $t_n^S(x)$ and (ii) $t_n^G(x)$, resp. As the sample size $n \to \infty$ we get the counterparts of the standard normal cdf $\Phi(x)$: (i) $\Phi^S(x) := \lim_{n\to\infty} t_n^S(x)$, and (ii) $\Phi^G(x) := \lim_{n\to\infty} t_n^G(x)$. Explicit formulae are given for the underlying new cdf's. For example $\Phi^G(x) = \Phi(x)$ iff $|x| \geq \sqrt{3}$. Thus the classical 95% confidence interval for the unknown expected value of Gaussian distributions covers the center of symmetry with at least 95% probability for Gaussian scale mixture distributions. On the other hand, the 90% quantile of $\Phi^G$ is $4\sqrt{3}/5 = 1.385\cdots > \Phi^{-1}(0.9) = 1.282\ldots$.

## 1. Introduction

An inspiring recent paper by Lehmann [9] summarizes Student's contributions to small sample theory in the period 1908–1933. Lehmann quoted Student [10]: "The question of applicability of normal theory to non-normal material is, however, of considerable importance and merits attention both from the mathematician and from those of us whose province it is to apply the results of his labours to practical work."

In this paper we consider two important classes of distributions. The first class is the class of all symmetric distributions. The second class consists of scale mixtures of normal distributions which contains all symmetric stable distributions, Laplace, logistic, exponential power, Student's $t$, etc. For scale mixtures of normal distributions see Kelker [8], Efron and Olshen [5], Gneiting [7], Benjamini [1]. Gaussian scale mixtures are important in finance, bioinformatics and in many other areas of applications where the errors are heavy tailed.

First, let $X_1, X_2, \ldots, X_n$ be independent (not necessarily identically distributed) observations, and let $\mu$ be an unknown parameter with

$$X_i = \mu + \xi_i, \ \ i = 1, 2, \ldots, n,$$

where the random errors $\xi_i$, $1 \leq i \leq n$ are independent, and symmetrically distributed around zero. Suppose that

$$\xi_i = s_i \eta_i, \quad i = 1, 2, \ldots, n,$$

where $s_i, \eta_i$ $i = 1, 2, \ldots, n$ are independent pairs of random variables, and the random scale, $s_i \geq 0$, is also independent of $\eta_i$. We also assume the $\eta_i$ variables are identically distributed with given cdf $F$ such that $F(x) + F(-x^-) = 1$ for all real numbers $x$.

[1]Department of Mathematics and Statistics, Bowling Green State University, Bowling Green, OH 43403-0221 and Alfréd Rényi Institute of Mathematics, Hungarian Academy of Sciences, Budapest, Hungary, e-mail: gabors@bgnet.bgsu.edu

Student's $t$-statistic is defined as $T_n = \sqrt{n}(\overline{X} - \mu)/S$, $n = 2, 3, \ldots$ where $\overline{X} = \sum_{i=1}^{n} X_i/n$ and $S^2 = \sum_{i=1}^{n}(X_i - \overline{X})^2/(n-1) \neq 0$.

Introduce the notation
$$a^2 := \frac{nx^2}{x^2 + n - 1}.$$

For $x \geq 0$,

$$P\{|T_n| > x\} = P\{T_n^2 > x^2\} = P\left\{\frac{\left(\sum_{i=1}^{n} \xi_i\right)^2}{\sum_{i=1}^{n} \xi_i^2} > a^2\right\}. \tag{1.1}$$

(For the idea of this equation see Efron [4] p. 1279.) Conditioning on the random scales $s_1, s_2, \ldots, s_n$, (1.1) becomes

$$P\{|T_n| > x\} = EP\left\{\frac{\left(\sum_{i=1}^{n} s_i \eta_i\right)^2}{\sum_{i=1}^{n} s_i^2 \eta_i^2} > a^2 | s_1, s_2, \ldots, s_n\right\}$$
$$\leq \sup_{\substack{\sigma_k \geq 0 \\ k=1,2,\ldots,n}} P\left\{\frac{\left(\sum_{i=1}^{n} \sigma_i \eta_i\right)^2}{\sum_{i=1}^{n} \sigma_i^2 \eta_i^2} > a^2\right\},$$

where $\sigma_1, \sigma_2, \ldots, \sigma_n$ are arbitrary nonnegative, *non-random* numbers with $\sigma_i > 0$ for at least one $i = 1, 2, \ldots, n$.

For Gaussian errors $P\{|T_n| > x\} = P(|t_{n-1}| > x)$ where $t_{n-1}$ is a t-distributed random variable with $n-1$ degrees of freedom. The corresponding cdf is denoted by $t_{n-1}(x)$. Suppose $a \geq 0$. For scale mixtures of the cdf $F$ introduce

$$1 - t_{n-1}^{(F)}(a) := \frac{1}{2} \sup_{\substack{\sigma_k \geq 0 \\ k=1,2,\ldots,n}} P\left\{\frac{\left(\sum_{i=1}^{n} \sigma_i \eta_i\right)^2}{\sum_{i=1}^{n} \sigma_i^2 \eta_i^2} > a^2\right\}. \tag{1.2}$$

For $a < 0$, $t_{n-1}^{(F)}(a) := 1 - t_{n-1}^{(F)}(-a)$. It is clear that if $1 - t_{n-1}^{(F)}(a) \leq \alpha/2$, then $P\{|T_n| > x\} \leq \alpha$. This is the starting point of our two excursions.

First, we assume $F$ is the cdf of a symmetric Bernoulli random variable supported on $\pm 1$ ($p = 1/2$). In this case the set of scale mixtures of $F$ is the complete set of symmetric distributions around 0, and the corresponding $t$ is denoted by $t^S$ ($t_n^S(x) = t_n^{(F)}(x)$ when $F$ is the Bernoulli cdf). In the second excursion we assume $F$ is Gaussian, the corresponding $t$ is denoted by $t^G$.

How to choose between these two models? If the error tails are lighter than the Gaussian tails, then of course we cannot apply the Gaussian scale mixture model. On the other hand, there are lots of models (for example the variance gamma model in finance) where the error distributions are supposed to be scale mixtures of Gaussian distributions (centered at 0). In this case it is preferable to apply the second model because the corresponding upper quantiles are smaller. For an intermediate model where the errors are symmetric and unimodal see Székely and Bakirov [11]. Here we could apply a classical theorem of Khinchin (see Feller [6]); according to this theorem all symmetric unimodal distributions are scale mixtures of symmetric uniform distributions.

## 2. Symmetric errors: scale mixtures of coin flipping variables

Introduce the Bernoulli random variables $\varepsilon_i$, $P(\varepsilon_i = \pm 1) = 1/2$. Let $\mathcal{P}$ denote the set of vectors $\mathbf{p} = (p_1, p_2, \ldots, p_n)$ with Euclidean norm 1, $\sum_{k=1}^n p_k^2 = 1$. Then, according to (1.2), if the role of $\eta_i$ is played by $\varepsilon_i$ with the property that $\varepsilon_i^2 = 1$,

$$1 - t_{n-1}^S(a) = \sup_{\mathbf{p} \in \mathcal{P}} P\{p_1\varepsilon_1 + p_2\varepsilon_2 + \cdots + p_n\varepsilon_n \geq a\}.$$

The main result of this section is the following.

**Theorem 2.1.** *For $0 < a \leq \sqrt{n}$, $2^{-\lceil a^2 \rceil} \leq 1 - t_{n-1}^S(a) = \frac{m}{2^n}$,*

*where $m$ is the maximum number of vertices $\mathbf{v} = \overbrace{(\pm 1, \pm 1, \ldots, \pm 1)}^{n}$ of the $n$-dimensional standard cube that can be covered by an $n$-dimensional closed sphere of radius $r = \sqrt{n - a^2}$. (For $a > \sqrt{n}$, $1 - t_{n-1}^S(a) = 0$.)*

*Proof.* Denote by $\mathcal{P}_a$ the set of all $n$-dimensional vectors with Euclidean norm $a$. The crucial observation is the following. For all $a > 0$,

$$\begin{aligned}
1 - t_{n-1}^S(a) &= \sup_{\mathbf{p} \in \mathcal{P}} P\left\{\sum_{j=1}^n p_j \varepsilon_j \geq a\right\} \\
&= \sup_{\mathbf{p} \in \mathcal{P}_a} P\left\{\sum_{j=1}^n (\varepsilon_j - p_j)^2 \leq n - a^2\right\}.
\end{aligned} \tag{2.1}$$

Here the inequality $\sum_{j=1}^n (\varepsilon_j - p_j)^2 \leq n - a^2$ means that the point

$$\mathbf{v} = (\varepsilon_1, \varepsilon_2, \ldots, \varepsilon_n),$$

a vertex of the $n$ dimensional standard cube, falls inside the (closed) sphere $G(\mathbf{p}, r)$ with center $\mathbf{p} \in \mathcal{P}_a$ and radius $r = \sqrt{n - a^2}$. Thus

$$1 - t_n^S(a) = \frac{m}{2^n},$$

where $m$ is the maximal number of vertices $\mathbf{v} = \overbrace{(\pm 1, \pm 1, \ldots, \pm 1)}^{n}$ which can be covered by an $n$-dimensional closed sphere with given radius $r = \sqrt{n - a^2}$ and varying center $\mathbf{p} \in \mathcal{P}_a$. It is clear that without loss of generality we can assume that the Euclidean norm of the optimal center is $a$.

If $k \geq 0$ is an integer and $a^2 \leq n - k$, then $m \geq 2^k$ because one can always find $2^k$ vertices which can be covered by a sphere of radius $\sqrt{k}$. Take, e.g., the vertices

$$(\overbrace{1,1,1,\ldots,1}^{n-k}, \overbrace{\pm 1, \pm 1, \ldots, \pm 1}^{k}),$$

and the sphere $G(\mathbf{c}, \sqrt{k})$ with center

$$\mathbf{c} = (\overbrace{1,1,\ldots,1}^{n-k}, \overbrace{0,0,\ldots,0}^{k}).$$

With a suitable constant $0 < C \leq 1$, $\mathbf{p} = C\mathbf{c}$ has norm $a$ and since the squared distances of $\mathbf{p}$ and the vertices above are $kC \leq k$, the sphere $G(\mathbf{p}, \sqrt{k})$ covers $2^k$ vertices. This proves the lower bound $2^{-\lceil a^2 \rceil} \leq 1 - t_n^S(a)$ in the Theorem. Thus the theorem is proved. □

**Remark 1.** Critical values for the $t^S$-test can be computed as the infima of the $x$-values for which $t^S_{n-1}\left(\sqrt{\frac{nx^2}{n-1+x^2}}\right) \leq \alpha$.

**Remark 2.** Define the counterpart of the standard normal distribution as follows.

$$\Phi^S(a) \stackrel{def}{=} \lim_{n \to \infty} t^S_n(a).$$

Theorem 1 implies that for $a > 0$,

$$1 - 2^{-\lceil a^2 \rceil} \leq \Phi^S(a). \tag{2.2}$$

Our computations suggest that the upper tail probabilities of $\Phi^S$ can be approximated by $2^{-\lceil a^2 \rceil}$ so well that the .9, .95, .975 quantiles of $\Phi^S$ are equal to $\sqrt{3}$, 2, $\sqrt{5}$, resp. with at least three decimal precision. We conjecture that $\Phi^S(\sqrt{3}) = .9$, $\Phi^S(2) = .95$, $\Phi^S(\sqrt{5}) = .975$. On the other hand, the .999 and higher quantiles almost coincide with the corresponding standard normal quantiles, thus in this case we do not need to pay a heavy price for dropping the condition of normality. On this problem see also the related papers by Eaton [2] and Edelman [3].

## 3. Gaussian scale mixture errors

An important subclass of symmetric distributions consists of the scale mixture of Gaussian distributions. In this case the errors can be represented in the form $\xi_j = s_i Z_i$ where $s_i \geq 0$ as before and independent of the standard normal $Z_i$. We have the equation

$$1 - t^G_{n-1}(a) = \sup_{\substack{\sigma_k \geq 0 \\ k=1,2,\ldots,n}} P\left\{\frac{\sigma_1 Z_1 + \sigma_2 Z_2 + \cdots + \sigma_n Z_n}{\sqrt{\sigma_1^2 Z_1^2 + \sigma_2^2 Z_2^2 + \cdots + \sigma_n^2 Z_n^2}} \geq a\right\}. \tag{3.1}$$

Recall that $a^2 = \frac{nx^2}{n-1+x^2}$ and thus $x = \sqrt{\frac{a^2(n-1)}{n-a^2}}$.

**Theorem 3.1.** *Suppose $n > 1$. Then for $0 \leq a < 1$, $t^G_{n-1}(a) = 1/2$, $t^G_{n-1}(1) = 3/4$, for $a \geq \sqrt{n}$, $t^G_{n-1}(a) = 1$, and finally, for $1 < a < \sqrt{n}$,*

$$1 - t^G_{n-1}(a) = \max_{1 \leq k \leq n} P\left\{\frac{Z_1 + Z_2 + \cdots + Z_k}{\sqrt{Z_1^2 + Z_2^2 + \cdots + Z_k^2}} \geq a\right\}$$

$$= \max_{a^2 < k \leq n} P\left(t_{k-1} > \sqrt{\frac{a^2(k-1)}{k-a^2}}\right)$$

*where $t_{k-1}$ is a t-distributed random variable with $k - 1$ degrees of freedom.*

The point of this theorem is that $\sup_{\sigma_1,\sigma_2,\ldots,\sigma_n}$ in (3.1) is taken when all nonzero $\sigma's$ are equal and here the number of zeros depends on $a$. For details see Székely and Bakirov [11].

Compute the intersection points of the curves

$$P\left(t_{k-1} > \sqrt{\frac{a^2(k-1)}{k-a^2}}\right)$$

for two neighboring indices. We get the following equation

$$\frac{2\Gamma\left(\frac{k}{2}\right)}{\sqrt{\pi(k-1)}\Gamma\left(\frac{k-1}{2}\right)} \int_0^{\sqrt{\frac{a^2(k-1)}{k-a^2}}} \left(1+\frac{u^2}{k-1}\right)^{-\frac{k}{2}} du$$
$$= \frac{2\Gamma\left(\frac{k+1}{2}\right)}{\sqrt{\pi k}\Gamma\left(\frac{k}{2}\right)} \int_0^{\sqrt{\frac{a^2 k}{k+1-a^2}}} \left(1+\frac{u^2}{k}\right)^{-\frac{k+1}{2}} du.$$

for the intersection point $A(k)$. It is not hard to show that $\lim_{k\to\infty} A(k) = \sqrt{3}$. This leads to the following:

**Corollary 1.** *There exists a sequence $A(1) := 1 < A(2) < A(3) < \cdots < A(k) \longrightarrow \sqrt{3}$, such that*

(i) *for $a \in [A(k-1), A(k)]$, $k = 2, 3, \ldots, n-1$,*

$$t_{n-1}^G(a) = P\left(t_{k-1} > \sqrt{\frac{a^2(k-1)}{k-a^2}}\right),$$

(ii) *for $a \geq \sqrt{3}$ that is for $x > \sqrt{3(n-1)/(n-3)}$,*

$$t_{n-1}^G(a) = t_{n-1}(a).$$

The most surprising part of Corollary 1 is of course the nice limit, $\sqrt{3}$. This shows that above $\sqrt{3}$ the usual $t$-test applies even if the errors are not necessarily normals only scale mixtures of normals. Below $\sqrt{3}$, however, the 'robustness' of the $t$-test gradually decreases. Splus can easily compute that $A(2) = 1.726$, $A(3) = 2.040$. According to our Table 1, the one sided 0.025 level critical values coincide with the classical t-critical values.

Recall that for $x \geq 0$, the Gaussian scale mixture counterpart of the standard normal cdf is

$$\Phi^G(x) := \lim_{n\to\infty} t_n^G(x) \tag{3.2}$$

(Note that in the limit, as $n \to \infty$, we have $a = x$ if both are assumed to be nonnegative; $\Phi^G(-x) = 1 - \Phi^G(x)$.)

**Corollary 2.** *For $0 \leq x < 1$, $\Phi^G(x) = .5$, $\Phi^G(1) = .75$, and for $x \geq \sqrt{3}$, $\Phi^G(x) = \Phi(x)$, where $\Phi(x)$ is the standard normal cdf ($\Phi^G(\sqrt{3}) = \Phi(\sqrt{3}) = 0.958$).*

For quantiles between .5 and .875 the max in Theorem 3.1 is taken at $k = 2$ and thus in this interval $\Phi^G(x) = C(x/\sqrt{(2-x^2)})$, where $C(x)$ is the standard Cauchy cdf. This is the convex section of the curve $\Phi^G(x), x \geq 0$. Interestingly the convex part is followed by a linear section: $\Phi^G(x) = x/(2\sqrt{3}) + 1/2$ for $1.3136 \cdots < x < 1.4282 \ldots$. Thus the 90% quantile is exactly $4\sqrt{3}/5$: $\Phi^G(4\sqrt{3}/5) = 0.9$. The following critical values are important in applications: $0.95 = \Phi(1.645) = \Phi^G(1.650)$, $0.9 = \Phi(1.282) = \Phi^G(1.386)$, $0.875 = \Phi(1.150) = \Phi^G(1.307)$ (see the last row of Table 1).

**Remark 3.** It is clear that for $a > 0$ we have the inequalities $t_n(a) \geq t_n^G(a) \geq t_n^S(a)$. According to Corollary 1, the first inequality becomes an equality iff $a \geq \sqrt{3}$. In connection with the second inequality one can show that the difference of the $\alpha$-quantiles of $t_n^G(a)$ and $t_n^S(a)$ tends to 0 as $\alpha \to 1$.

TABLE 1
*Critical values for Gaussian scale mixture errors computed from $t_n^G(\sqrt{nx^2/(n-1+x^2)}) = \alpha$*

| $n-1$ | 0.125 | 0.100 | 0.050 | 0.025 |
|---:|---:|---:|---:|---:|
| 2 | 1.625 | 1.886 | 2.920 | 4.303 |
| 3 | 1.495 | 1.664 | 2.353 | 3.182 |
| 4 | 1.440 | 1.579 | 2.132 | 2.776 |
| 5 | 1.410 | 1.534 | 2.015 | 2.571 |
| 6 | 1.391 | 1.506 | 1.943 | 2.447 |
| 7 | 1.378 | 1.487 | 1.895 | 2.365 |
| 8 | 1.368 | 1.473 | 1.860 | 2.306 |
| 9 | 1.361 | 1.462 | 1.833 | 2.262 |
| 10 | 1.355 | 1.454 | 1.812 | 2.228 |
| 11 | 1.351 | 1.448 | 1.796 | 2.201 |
| 12 | 1.347 | 1.442 | 1.782 | 2.179 |
| 13 | 1.344 | 1.437 | 1.771 | 2.160 |
| 14 | 1.341 | 1.434 | 1.761 | 2.145 |
| 15 | 1.338 | 1.430 | 1.753 | 2.131 |
| 16 | 1.336 | 1.427 | 1.746 | 2.120 |
| 17 | 1.335 | 1.425 | 1.740 | 2.110 |
| 18 | 1.333 | 1.422 | 1.735 | 2.101 |
| 19 | 1.332 | 1.420 | 1.730 | 2.093 |
| 20 | 1.330 | 1.419 | 1.725 | 2.086 |
| 21 | 1.329 | 1.417 | 1.722 | 2.080 |
| 22 | 1.328 | 1.416 | 1.718 | 2.074 |
| 23 | 1.327 | 1.414 | 1.715 | 2.069 |
| 24 | 1.326 | 1.413 | 1.712 | 2.064 |
| 25 | 1.325 | 1.412 | 1.709 | 2.060 |
| 100 | 1.311 | 1.392 | 1.664 | 1.984 |
| 500 | 1.307 | 1.387 | 1.652 | 1.965 |
| 1,000 | 1.307 | 1.386 | 1.651 | 1.962 |

Our approach can also be applied for two-sample tests. In a joint forthcoming paper with N. K. Bakirov the Behrens–Fisher problem will be discussed for Gaussian scale mixture errors with the help of our $t_n^G(x)$ function.

## Acknowledgments

The author also wants to thank many helpful suggestions of N. K. Bakirov, M. Rizzo, the referees of the paper, and the editor of the volume.